\documentclass[preprint,12pt]{elsarticle}

\usepackage{amssymb}
\usepackage{amsmath}
\usepackage{amsthm}
\usepackage[T1]{fontenc}
\usepackage[french,english]{babel}
\usepackage[colorlinks]{hyperref}
\AtBeginDocument{%
  \hypersetup{%
    linkcolor=blue,%
    citecolor=green,%
  }%
}
\usepackage{cleveref}
\usepackage{geometry}

\crefname{equation}{}{}

\newtheorem{theorem}{Theorem}[section]
\newtheorem{lemma}[theorem]{Lemma}

\newtheorem{corollary}[theorem]{Corollary}

\theoremstyle{definition}
\newtheorem{definition}[theorem]{Definition}

\newtheorem{condition}[theorem]{Condition}

\theoremstyle{remark}
\newtheorem{remark}[theorem]{Remark}

\numberwithin{equation}{section}

\journal{Journal de Math\'{e}matiques Pures et Appliqu\'{e}es}

\newenvironment{abstracts}
 {\global\setbox\absbox=\vbox\bgroup
    \hsize=\textwidth
    \linespread{1}\selectfont}
 {\vspace{-\bigskipamount}\egroup}
\renewenvironment{abstract}[1][]
 {\if\relax\detokenize{#1}\relax\else\selectlanguage{#1}\fi
  \noindent\textbf{\abstractname}\par\medskip\noindent\ignorespaces}
 {\par\bigskip}

\begin{document}

\begin{frontmatter}



\title{Boundary H\"{o}lder Regularity for Elliptic Equations\tnoteref{t1}}

\author[rvt2]{Yuanyuan Lian}
\ead{lianyuanyuan@nwpu.edu.cn; lianyuanyuan.hthk@gmail.com}
\author[rvt2]{Kai Zhang\corref{cor1}}
\ead{zhang\_kai@nwpu.edu.cn; zhangkaizfz@gmail.com}
\author[rvt]{Dongsheng Li}
\ead{lidsh@mail.xjtu.edu.cn}
\author[rvt]{Guanghao Hong}
\ead{ghhongmath@mail.xjtu.edu.cn}

\tnotetext[t1]{This research is supported by the National Natural Science Foundation of China (Grant No. 11701454 and 11671316) and the Natural Science Basic Research Plan in Shaanxi Province of China (Program No. 2018JQ1039).}

\cortext[cor1]{Corresponding author. ORCID: \href{https://orcid.org/0000-0002-1896-3206}{0000-0002-1896-3206}}

\address[rvt2]{School of Mathematics and Statistics, Northwestern Polytechnical University, Xi'an, Shaanxi, 710129, PR China}
\address[rvt]{School of Mathematics and Statistics, Xi'an Jiaotong University, Xi'an, Shaanxi, 710049, PR China}

\begin{abstracts}

\begin{abstract}
This paper investigates the relation between the boundary geometric properties and the boundary regularity of the solutions of elliptic equations. We prove by a new unified method the pointwise boundary H\"{o}lder regularity under proper geometric conditions. ``Unified'' means that our method is applicable for the Laplace equation, linear elliptic equations in divergence and non-divergence form, fully nonlinear elliptic equations, the $p-$Laplace equations and the fractional Laplace equations etc. In addition, these geometric conditions are quite general. In particular, for local equations, the measure of the complement of the domain near the boundary point concerned could be zero. The key observation in the method is that the strong maximum principle implies a decay for the solution, then a scaling argument leads to the H\"{o}lder regularity. Moreover, we also give a geometric condition, which guarantees the solvability of the Dirichlet problem for the Laplace equation. The geometric meaning of this condition is more apparent than that of the Wiener criterion.
\end{abstract}

\begin{abstract}[french]
Dans cet article, nous \'{e}tudies la relation entre les propri\'{e}t\'{e}s g\'{e}om\'{e}triques des fronti\`{e}re et la r\'{e}gularit\'{e} fronti\`{e}re des solutions d'\'{e}quations elliptiques. Nous prouvons par une nouvelle m\'{e}thode unifi\'{e}e la r\'{e}gularit\'{e} h\"{o}ld\'{e}rienne ponctuelle dans des conditions g\'{e}om\'{e}triques appropri\'{e}es. Unifi\'{e} signifie que notre m\'{e}thode est applicable \`{a} l'\'{e}quation de Laplace, aux \'{e}quations elliptiques lin\'{e}aires sous forme de divergence et de non-divergence, aux \'{e}quations elliptiques enti\`{e}rement non lin\'{e}aires, aux \'{e}quations p-Laplace et aux \'{e}quations fractionnelles de Laplace, etc. En outre, ces conditions g\'{e}om\'{e}triques sont assez g\'{e}n\'{e}rales. En particulier, pour les \'{e}quations locales, la mesure du compl\'{e}ment du domaine pr\`{e}s du point fronti\`{e}re concern\'{e} pourrait \^{e}tre nulle. Une observation cl\'{e} de notre m\'{e}thode est que le principe du maximum fort implique une d\'{e}croissance de la solution, puis un argument d'\'{e}chelle nous conduit \`{a} la r\'{e}gularit\'{e} de H\"{o}lder. De plus, nous donnons \'{e}galement une condition g\'{e}om\'{e}trique, qui garantit la solvabilit\'{e} du probl\`{e}me de Dirichlet pour l'\'{e}quation de Laplace. La signification g\'{e}om\'{e}trique de cette condition est plus apparente que celle du crit\`{e}re de Wiener.
\end{abstract}

\end{abstracts}

\end{frontmatter}


\section{Introduction}\label{S1}
Let $\Omega\subset R^{n}$ be a bounded domain and $g\in C(\partial \Omega) $. It has been taken for granted for a time that there always exists $u\in C(\bar{\Omega})$ such that $u$ is harmonic in $\Omega$ and $u\equiv g$ on $\partial \Omega$. That is, the Dirichlet problem for the Laplace equation is solvable for any bounded domain. However, in 1913, Lebesgue \cite{Lebesgue1912} constructed a bounded domain on which the Dirichlet problem is not solvable. This indicates that the domain must satisfy some condition for the continuity of the solution up to the boundary. In 1924, Wiener \cite{Wiener-1924} proposed a sufficient and necessary condition for the solvability of the Dirichlet problem. This is the famous Wiener criterion which solves the Dirichlet problem completely. The Wiener criterion has been extended to the linear equations in divergence form \cite{L-S-W} and quasilinear equations including the $p$-Laplace equations \cite{L-M,K-M,Mazja}.

However, there are some disadvantages in the Wiener criterion. It is not easy to check whether a domain satisfies the Wiener criterion. The notion capacity is used and calculating the capacity of a set is difficult in general. In addition, the generalization to other types of equations is limited because the definition of capacity is close to the divergence structure of the equation. Moreover, there is no continuity modulus estimate in the Wiener criterion. It doesn't point out which kind of continuity up to the boundary for the solution. Quantitative estimates for continuity modulus are important in the regularity theory.

The H\"{o}lder continuity is a kind of quantitative estimate. It is usually the first smooth regularity for solutions and the beginning for higher regularity. It can also provide compactness in some problems. With respect to the boundary H\"{o}lder continuity, we mention the following results in which geometric conditions are given and quantitative estimates are also derived. If $\Omega$ satisfies the exterior cone condition at $x_0\in \partial \Omega $, then the solution is H\"{o}lder continuous at $x_0$ (see \cite[Problem 2.12]{MR1814364} and \cite{MR0221087}). If $\Omega$ satisfies the exterior sphere condition at $x_0$, then the solution is Lipschitz continuous at $x_0$. The later one has been generalized to exterior Dini hypersurface condition (see \cite{MR3167627} and \cite{Safonov_2008}). To the best of our knowledge, the exterior cone condition is the weakest geometric condition for boundary H\"{o}lder regularity.

In this paper, we provide a new method to prove the boundary H\"{o}lder regularity. This method is not only appropriate for the Laplace equation but also applicable for other kinds of equations, including linear elliptic equations in divergence form and non-divergence form, fully nonlinear elliptic equations, the $p-$Laplace equations and the fractional Laplace equations etc. We will propose several conditions with clear geometric meaning and then prove the pointwise boundary H\"{o}lder regularity for the corresponding equations. These geometric conditions are generalized widely from the exterior cone condition. In particular, for local equations, the measure of the complement of the domain near a boundary point could be zero. Finally, we also give a geometric condition for the Laplace equation to guarantee the solvability of the Dirichlet problem. We remark here that the equations considered in this paper are only some concrete examples. This method may have a wide range of potential applications.

Now, we clarify the key idea briefly. Instead of proving the boundary H\"{o}lder regularity by constructing a (local) barrier or applying the Harnack inequality at the boundary, we solve the Dirichlet problems in a sequence of balls centered at the boundary point concerned. Then by applying the strong maximum principle, the comparison principle and the scaling invariant property of the equations, we derive a quantitative decay of the oscillation of the solution near the boundary point, which implies the boundary H\"{o}lder continuity immediately. These properties used above are occupied by many types of equations. Hence, this method is easily extended to other types of equations.

This paper is organized as  follows. In Section 2, we present the geometric conditions on the domains. In Section 3, we prove the pointwise boundary H\"{o}lder regularity for different elliptic equations under the corresponding conditions proposed in Section 2. The solvability of the Dirichlet problem for the Laplace equation will be proved in the last section.

\section{Geometric conditions}
As is well known, the geometric property of the domain near some boundary point has significant influence on the boundary regularity there. This is one of the most important difference from the interior regularity. In this section, we introduce some geometric conditions on the domains. The corresponding boundary regularity will be proved in later sections.
\begin{definition}\label{d-1.1}
Let $\{r_k\}_{k=0}^{\infty}$ be  a positive sequence. We call it a quasi-geometric sequence if there exist constants $0<\tau_1<\tau_2<1$ such that
\begin{equation}\label{e-rk-1}
    \tau_1 r_{k-1}\leq r_k \leq \tau_2 r_{k-1}, ~~\forall ~ k\geq 1.
\end{equation}
\end{definition}
\begin{remark}
If there exist $0<\tau<1$ and a positive sequence $\{r_k\}_{k=0}^{\infty}$ such that $r_k\rightarrow 0$ and $\tau r_{k-1}\leq r_k$ ($k\geq 1$), then it is easy to verify that there exists a subsequence of $\{r_k\}$ satisfying\cref{e-rk-1} with $\tau_1=\tau^2$ and $\tau_2=\tau$. Hence, the essence of\cref{e-rk-1} is that $r_k$ should not decrease too rapidly.
\end{remark}

The following geometric condition will be used to prove the boundary H\"{o}lder regularity for the Laplace equation.
\begin{condition}[\textbf{H1}]\label{d-gc-1}
Let $\Omega\subset R^{n}$ be a bounded domain and $x_0\in \partial \Omega$. We say that $\Omega$ satisfies the (H1) condition at $x_0$ if there exist a constant $0<\nu<1$ and a quasi-geometric sequence $\{r_k\}_{k=0}^{\infty}$ such that
\begin{equation}\label{e-h2-2-2}
 \frac{H^{n-1}( \partial B(x_0,r_k) \cap \Omega ^c)}{r_k^{n-1}}\geq \nu,~~\forall ~ k\geq 0,
\end{equation}
where $H^{n-1}$ denotes the $n-1$ dimensional Hausdorff measure.
\end{condition}

For general equations, there is no explicit expression connecting the solutions with the boundary values. Nevertheless, the strong maximum principle and the scaling property implies a quantitative decay for the solutions. Hence, we can obtain the boundary H\"{o}lder regularity for general equations. The following stronger (compared to the (H1) condition) geometric condition will be used to prove boundary H\"{o}lder regularity for general (local) elliptic equations.
\begin{condition}[\textbf{H2}]\label{d-gc}
Let $\Omega\subset R^{n}$ be a bounded domain and $x_0\in \partial \Omega$. We say that $\Omega$ satisfies the (H2) condition at $x_0$ if there exist a constant $0<\nu<1$, a quasi-geometric sequence $\{r_k\}_{k=0}^{\infty}$ and a sequence $\{y_k\}_{k=0}^{\infty}$ with $y_k\in \partial B(x_0,r_k)$ such that
\begin{equation}\label{e-bnu}
  \partial B(x_0,r_k) \cap B(y_k,\nu r_k) \subset \Omega ^c,~~\forall ~k\geq 0.
\end{equation}
\end{condition}

Before proceeding to present other geometric conditions, we make some remarks. We say that $\Omega$ satisfies the exterior cone condition at $x_0\in \partial \Omega$ if there exists a finite right circular cone $K$ with vertex $x_0$ such that $\bar{\Omega}\cap \bar{K}=x_0$. And, it is well known that the exterior cone condition is sufficient for the boundary H\"{o}lder regularity for the uniformly elliptic equations (see \cite[Problem 2.12]{MR1814364} and \cite{MR0221087}). Clearly, a domain satisfying the exterior cone condition at some boundary point satisfies the (H2) condition at the same point with $\tau_1=\tau_2=1/2$ and $\nu$ depending on the aperture of the cone. Hence, the (H2) condition is a generalization of the exterior cone condition.

In addition, a Reifenberg flat domain also satisfies the (H2) condition. We introduce the following definition for the Reifenberg flat domain (see \cite[Definition 2.5]{Byun-Wang-2008}).
\begin{definition}[\textbf{Reifenberg flat domain}]\label{d-re}
We say that $\Omega$ satisfies the exterior ($\delta$,R)-Reifenberg ($\delta<1$) flat condition at $x_0\in \partial \Omega$ if for any $0<r<R$, there exists a coordinate system $\{y_1,...,y_n \}$ such that $x_0=0$ in this coordinate system and
\begin{equation}\label{e-re}
B_r\cap \Omega \subset B_r \cap\{y_n >-\delta r\}.
\end{equation}
\end{definition}

Suppose that $\Omega$ satisfies the exterior ($\delta$,R)-Reifenberg flat condition at some boundary point $x_0$. By\cref{e-re}, for any $0<r<R$, in the new coordinate, $\partial B_r \cap\{y_n \leq-\delta r\}\subset \Omega^c$. Hence, $\Omega$ satisfies the (H2) condition at $x_0$ with $\tau_1=\tau_2=1/2$ and $\nu$ depending on $\delta$. The Reifenberg flat domain was first introduced by Reifenberg in 1960 \cite{MR0114145} and appears in minimal surface theory and free boundary problems. One interesting feature of a Reifenberg flat domain is that its boundary could be a fractal. For more research on Reifenberg flat domains, we refer to \cite{MR2069724, Byun-Wang-2008} and the references therein.

More generally, the corkscrew domains, including the non-tangentially accessible domains (NTA domains for short) as a subclass, also satisfy the (H2) condition. The following is the definition for the corkscrew domains (see \cite[(3.1)]{MR676988}).
\begin{definition}[\textbf{Corkscrew domain}]\label{d-cork}
We say that $\Omega$ satisfies the exterior corkscrew condition at $x_0\in \partial \Omega$ if there exist $0<\delta<1/4$ and $R>0$ such that for any $0<r<R$, there exists $y\in B(x_0,r)$ such that $B(y, \delta r)\subset \Omega ^c$.
\end{definition}

Suppose that $\Omega$ satisfies the exterior corkscrew condition at $x_0\in \partial \Omega$. Then for any $0<r<R$, there exists $y\in B(x_0,r)$ such that $B(y, \delta r)\subset \Omega ^c$. Let $r_1=|y-x_0|$. Then $\delta r\leq r_1\leq (1-\delta) r$ and $y\in \partial B(x_0,r_1)$. Moreover, $\partial B(x_0,r_1)\cap B(y, \delta r_1)\subset \partial B(x_0,r_1)\cap B(y, \delta r)\subset\Omega ^c$. Hence, $\Omega$ satisfies the (H2) condition at $x_0$ with $\tau_1=\delta$, $\tau_2=1-\delta$ and $\nu=\delta$. We remark here that the NTA domains and the corkscrew domains appear frequently in the study of harmonic measures and non-tangential limits (see \cite{MR676988,MR1446617}).

In fact, in the (H2) condition, we only need a sequence of uniform portion of the $(n-1)$ dimensional spheres contained in $\Omega^c$ rather than an $n$ dimensional set. In particular, the measure of the complement of the domain near the boundary point could be zero. For example, let $\Gamma=\partial B_{1/2}\cap B(y_0,\nu)$ where $y_0\in \partial B_{1/2}$ and $0<\nu<1$. Set
\begin{equation*}
  \Omega=B_1\backslash \{0\}\backslash \bigcup_{k=0}^{\infty} \frac{1}{2^k}\Gamma.
\end{equation*}
Then $\Omega$ satisfies the (H2) condition at $0\in \partial \Omega$ and the H\"{o}lder regularity at $0$ for elliptic equations can be proved (see Section 2). Note that the Lebesgue measure of $B_{1/2}\cap \Omega^c$ is zero.

In this paper, we also consider the fractional Laplace equations. Since they are nonlocal equations, the condition (H2) is not appropriate here. Instead, we give the following geometric condition.
\begin{condition}[\textbf{H3}]\label{con-frac}
Let  $\Omega\subset R^{n}$ be a bounded domain and $x_0\in \partial \Omega$. We say that $\Omega$ satisfies the (H3) condition at $x_0$ if there exist a constant $0<\nu<1$ and a quasi-geometric sequence $\{r_k\}_{k=0}^{\infty}$ such that
\begin{equation}\label{e-nu-frac}
\frac{ |(B(x_0,r_k)\backslash B(x_0,r_{k+1}))\cap\Omega ^c|}{r_k^{n}}\geq \nu,~~\forall ~ k\geq 0,
\end{equation}
where $|\cdot|$ denotes the $n$ dimensional Lebesgue measure.
\end{condition}

\begin{remark}\label{r-con-frac}
It is also easy to verify that if a domain satisfies the exterior cone condition, the exterior $(\delta,R)$-Reifenberg flat condition, or the exterior corkscrew condition, it satisfies the (H3) condition. Hence, there is a large class of domains on which the boundary H\"{o}lder regularity holds for the fractional Laplace equations. Here, the measure of the complement of the domain near the boundary point can not be zero, which is a difference from the local equations.
\end{remark}

Finally, for the Laplace equation, there exists an explicit relationship between the solution and the boundary values via the Poisson integral. Hence, we have clear quantitative estimate for the decay of the solution around some boundary point. This allows us to obtain other boundary continuity under the corresponding geometric conditions. We give the following geometric condition for example and prove the boundary continuity of the solution in the last section.

\begin{condition}[\textbf{H4}]\label{con-4}
Let $\Omega\subset R^{n}$ be a bounded domain and $x_0\in \partial \Omega$. We say that $\Omega$ satisfies the (H4) condition at $x_0$ if there exist a quasi-geometric sequence $\{r_k\}_{k=0}^{\infty}$ such that
\begin{equation}\label{e-h4-2}
 \sum_{k=0}^{\infty}\frac{H^{n-1}( \partial B(x_0,r_k) \cap \Omega ^c)}{r_k^{n-1}}=+\infty.
\end{equation}
\end{condition}

\begin{remark}
We will prove the continuity up to $x_0$ under\cref{e-h4-2}. The famous Wiener criterion states (see \cite[Chap. 2.9]{MR1814364} and \cite{Wiener-1924}) that for the Laplace equation, the solution continues up to $x_0$ if and only if
\begin{equation}\label{e-wiener}
   \sum_{k=0}^{\infty}\frac{\mathrm{cap} \left(B(x_0,r^{k})\cap\Omega^c\right)}{r^{k(n-2)}}=+\infty,
\end{equation}
where $0<r<1$ and $\mathrm{cap}$ denotes the capacity which is defined as the following (see \cite[Chap. 2.9]{MR1814364}). For any $\Omega\subset R^n$, define
\begin{equation*}
  \mathrm{cap} ~\Omega=\underset{v\in K}{\inf} \int |Dv|^2,
\end{equation*}
where
\begin{equation*}
  K=\left\{v\in C^{1}_0 (R^n)|v=1 \mbox{~on~} \Omega\right\}.
\end{equation*}

 Note that the form of\cref{e-h4-2} is similar to that of\cref{e-wiener}. However,\cref{e-h4-2} shows apparent geometric meaning whereas the capacity is used in\cref{e-wiener}.
\end{remark}

In the rest of this paper, if we say that a domain satisfies the (H1), (H2), (H3) or (H4) condition at some boundary point, it indicates the quasi-geometric sequence $\{r_k\}_{k=0}^{\infty}$ with $0<\tau_1<\tau_2<1$ and $r_0=1$, the constant $0<\nu<1$, and the sequence $\{y_k\}_{k=0}^{\infty}$.

\section{Boundary H\"{o}lder regularity}

Let $\Omega \subset R^{n}$ be a bounded domain and $f$ is a function defined on $\bar{\Omega}$. We say that $f$ is $C^{\alpha}$ at $x_0\in \bar{\Omega}$ or $f\in C^{\alpha}(x_0)$ if there exists a constant $C$ such that
\begin{equation*}
  |f(x)-f(x_0)|\leq C|x-x_0|^{\alpha},~~\forall~x\in \bar{\Omega}.
\end{equation*}
Then, define $[f]_{C^{\alpha}(x_0)}=C$ and $\|f\|_{C^{\alpha}(x_0)}=\|f\|_{L^{\infty}(\Omega)}+[f]_{C^{\alpha}(x_0)}$.

In the following, we will prove the boundary H\"{o}lder regularity for different elliptic equations under the corresponding geometric conditions presented in the last section. The idea and the sketch of the proofs have been motivated by \cite{Li-Zhang-Oblique}. For simplicity, we assume that all solutions concerned are continuous up to the boundary throughout this paper. In fact, this assumption is not essential and we may derive the existence of a solution which is H\"{o}lder continuous up the boundary point at which the corresponding geometric condition is satisfied (see \Cref{t-final}).

To demonstrate our idea clearly, we first consider the Laplace equation in a simple form. Although the problem is simple, the proof has contained the essential ingredients for the boundary H\"{o}lder regularity.
\begin{theorem}\label{t-Laplace}
Suppose that $\Omega$ satisfies the (H1) condition at $0\in \partial \Omega$. Let $0\leq u \leq 1$ satisfy
\begin{equation*}
\left\{\begin{aligned}
  \Delta u&=0 ~~\mbox{in}~~\Omega; \\
  u&=0 ~~\mbox{on}~~\partial \Omega\cap B_{1}.
\end{aligned}\right.
\end{equation*}
Then $u$ is $C^{\alpha}$ at $0$ and
\begin{equation*}
  u(x)\leq 2|x|^{\alpha}, ~~\forall ~x\in \Omega\cap B_1,
\end{equation*}
where $0<\alpha<1$ depends only on $n,\tau_1, \tau_2$ and $ \nu$. Here, $\tau_1, \tau_2$ and $\nu$ are constants from \Cref{d-1.1} and the (H1) condition.
\end{theorem}
\proof  Let
\begin{equation*}
 g(x) \equiv\left\{\begin{aligned}
&0 ~~&&\mbox{on}~~ \partial B_{1}\cap \Omega^c;\\
&1 ~~&&\mbox{on}~~\partial B_{1}\cap \Omega
\end{aligned}\right.
\end{equation*}
and $v$ be the Poisson integral of $g$ on $B_{1}$, i.e.,
\begin{equation*}
  v(x)=\frac{1-|x|^2}{n\omega _n} \int_{\partial B_{1}}\frac{g(y)ds}{|x-y|^n},
\end{equation*}
where $x\in B_1$ and $\omega_n$ denotes the volume of $B_1$.

Then $v$ is a positive harmonic function in $B_{1}$ and thus
\begin{equation*}
  \liminf_{x\rightarrow \partial \Omega\cap \bar{B}_1} v \geq 0.
\end{equation*}
Since $\Omega$ is open, $g$ is continuous at any point $x\in \partial B_1\cap \Omega$ ($g\equiv 1$ near $x$). Note that $v$ is the Poisson integral of $g$. Hence, $v$ is continuous up to $x$ and $v(x)=1$ (see the proof of \cite[Theorem 2.6]{MR1814364} for the continuity of $v$ at $x$). Then, we have
\begin{equation*}
  \liminf_{x\rightarrow \partial ( \Omega \cap B_1)} (v-u) \geq 0.
\end{equation*}
By the maximum principle (see \cite[Theorem 3.1]{MR1814364}),
\begin{equation}\label{e-Lap-comp-1}
 v-u\geq 0 ~~\mbox{on}~~\Omega\cap B_{1}.
\end{equation}
On the other hand, by the (H1) condition and the definition of $g$,
\begin{equation}\label{e-Lap-str}
    v\leq 1-\mu~~\mbox{in}~~ B_{\tau_2},
\end{equation}
where $0<\mu<1$ depends only on $n,\tau_2$ and $\nu$. Combining\cref{e-Lap-comp-1} and\cref{e-Lap-str}, we have
\begin{equation*}
  \sup_{\Omega\cap B_{r_1}} u\leq   \sup_{\Omega\cap B_{\tau_2}} v\leq 1-\mu.
\end{equation*}

Now, we prove the following by induction
\begin{equation}\label{e-Lap-decay-2}
  \sup_{\Omega\cap B_{r_k}} u \leq (1-\mu)^{k},~~\forall ~k\geq 1.
\end{equation}
In above argument, we have proved that it holds for $k=1$. Suppose that it holds for $k$ and we need to prove it for $k+1$. Let $y=x/r_k$, $w(y)=u(x)/(1-\mu)^{k}$ and $\tilde{\Omega}=\Omega/r_k$. Then $w$ and $\tilde{\Omega}$ satisfy the conditions of this theorem. Hence, from above argument, we have
\begin{equation*}
  \sup_{\Omega\cap B_{\tau_2}} w \leq 1-\mu.
\end{equation*}
By rescaling back to $u$, we have
\begin{equation*}
  \sup_{\Omega\cap B_{\tau_2r_k}} u \leq (1-\mu)^{k+1}.
\end{equation*}
Hence,
\begin{equation*}
  \sup_{\Omega\cap B_{r_{k+1}}} u \leq\sup_{\Omega\cap B_{\tau_2r_k}} u \leq (1-\mu)^{k+1}.
\end{equation*}

Then\cref{e-Lap-decay-2} implies the H\"{o}lder continuity of $u$ at $0$. Indeed, for any $x\in \Omega\cap B_{1}$, there exists $k$ such that $r_{k+1}\leq |x| \leq r_{k}$. Then by letting $1-\mu=\tau_1^{\alpha}$,
\begin{equation*}
u(x)\leq (1-\mu)^{k}= \frac{\tau_1^{(k+1) \alpha}}{1-\mu}\leq \frac{r_{k+1}^{\alpha}}{(1-\mu)}\leq
\frac{|x|^{\alpha}}{(1-\mu)}\leq 2|x|^{\alpha}.
\end{equation*}\qed

\begin{remark}
It is clear from the proof that larger $H^{n-1}(\partial B_1\cap\Omega^c)$ implies bigger $\mu$, which leads to faster oscillation decay for $u$. That is, $u$ has higher continuity. This gives the explanation that a more portion of $\Omega^c$ near a boundary point implies a higher continuity up to the boundary. For example,  the exterior cone condition implies the boundary H\"{o}lder continuity; the exterior sphere condition implies the boundary Lipschitz continuity.
\end{remark}

Now we prove the full result for the Laplace operator with the inhomogeneous boundary condition:
\begin{equation}\label{e-Poi}
\left\{\begin{aligned}
 &\Delta u=f &&\mbox{in}~~\Omega; \\
 & u=g &&\mbox{on}~~\partial\Omega.
\end{aligned}
\right.
\end{equation}
\begin{theorem}\label{t-Poi}
Suppose that $\Omega$ satisfies the (H1) condition at $0\in \partial \Omega$. Let $u$ satisfy\cref{e-Poi}
where $f\in L^p(\Omega)$ with $p>n/2$ and $g$ is $C^{\alpha}$ at $0$. Then $u$ is $C^{\beta}$ at $0$ and
\begin{equation}\label{h2}
  |u(x)-u(0)|\leq C |x|^{\beta}\left(\|u\|_{L^{\infty }(\Omega\cap B_1)}+\|f\|_{L^{p}(\Omega\cap B_1)}+[g]_{C^{\alpha}(0)}\right), ~~\forall ~x\in \Omega \cap B_1,
\end{equation}
where $0<\beta\leq \min(2-n/p,\alpha)$ depends only on $n,\tau_1, \tau_2$ and $\nu$; $C$ depends only on $n,\tau_2,\nu$ and $p$. Here, $\tau_1, \tau_2$ and $\nu$ are constants from \Cref{d-1.1} and the (H1) condition.
\end{theorem}
\proof Without loss of generality, we assume that $g(0)=0$. Otherwise, we may consider $u-g(0)$ instead. Let $M=\|u\|_{L^{\infty }(\Omega\cap B_1)}+\|f\|_{L^{p}(\Omega\cap B_1)}+[g]_{C^{\alpha}(0)}$ and $\Omega_{r}=\Omega\cap B_{r}$. To prove\cref{h2}, we only need to prove the following:

There exist constants $0<\beta\leq \min(2-n/p,\alpha)$ depending only on $n,\tau_1, \tau_2$ and $\nu$, and $\hat{C}$ depending only on $n,\tau_2,\nu$ and $p$ such that for all $k\geq 0$,

\begin{equation}\label{e-poisson-discrete}
\|u\|_{L^{\infty }(\Omega _{r_{k}})}\leq \hat{C} M r_{k}^{\beta}.
\end{equation}

We prove\cref{e-poisson-discrete} by induction. For $k=0$, it holds clearly. Suppose that it holds for $k$. We need to prove that it holds for $k+1$.

Take the zero extension of $f$ to the whole $R^{n}$ (similarly hereinafter) and let
\begin{equation*}
  \tilde{g}(x)\equiv \left\{\begin{aligned}
&Mr_{k}^{\alpha}~~&&\mbox{on}~~\partial B_{r_k}\cap \Omega^c;\\
&\hat{C} M r_{k}^{\beta}~~&&\mbox{on}~~\partial B_{r_k}\cap \Omega.
  \end{aligned}\right.
\end{equation*}
Define
\begin{equation*}
  v(x)=\frac{r_k^2-|x|^2}{n\omega _n r_k} \int_{\partial B_{r_k}}\frac{\tilde{g}(y)ds}{|x-y|^n} + \int_{B_{r_k}} G(x,y)(-|f|(y)) dy,
\end{equation*}
where $x\in B_{r_k}$ and $G$ is the Green's function for $B_{r_k}$. Then
\begin{equation*}
  \Delta (v-u)\leq 0 ~~\mbox{in}~~\Omega_{r_k}.
\end{equation*}

Next, similar to the proof of \Cref{t-Laplace}, we have
\begin{equation}\label{e.1}
  \liminf_{x\rightarrow \partial \Omega\cap \bar{B}_{r_k}} v \geq Mr_{k}^{\alpha},
\end{equation}
since $v\geq Mr_{k}^{\alpha}$ in $B_{r_k}$. Note also that $v$ is continuous and equals to $\hat{C} M r_{k}^{\beta}$ at any point $x\in \partial B_{r_k}\cap \Omega$. Hence, combining with\cref{e-poisson-discrete} and\cref{e.1}, we have
\begin{equation*}
  \liminf_{x\rightarrow \partial \Omega_{r_k}} (v-u) \geq 0.
\end{equation*}
By the maximum principle, we have
\begin{equation*}
u \leq v ~~\mbox{in}~~\Omega_{r_k}.
\end{equation*}

For $x\in B_{r_{k+1}}$, we have as before
\begin{equation}\label{e.2}
  \frac{r_k^2-|x|^2}{n\omega _n r_k} \int_{\partial B_{r_k}}\frac{\tilde{g}(y)ds}{|x-y|^n} \leq (1-\mu)\left(\hat{C}M r_{k}^{\beta}-M r_{k}^{\alpha}\right)+M r_{k}^{\alpha},
\end{equation}
where $0<\mu<1$ depends only on $n,\tau_2$ and $\nu$. It can also be verified easily that
\begin{equation*}
  \int_{B_{r_k}} |G(x,y)f(y)| dy\leq CMr_{k}^{2-n/p},
\end{equation*}
where $C$ depends only on $n$ and $p$. Hence,
\begin{equation*}
\begin{aligned}
  \|v\|_{L^{\infty}(B_{r_{k+1}})}&\leq (1-\mu)\left(\hat{C}M r_{k}^{\beta}-M r_{k}^{\alpha}\right)+M r_{k}^{\alpha}+CMr_{k}^{2-n/p}\\
&\leq \hat{C}M r_{k+1}^{\beta}\cdot \frac{1-\mu}{\tau_1^{\beta }}+\mu M r_{k}^{\beta}+ CMr_{k}^{\beta}\\
&\leq \hat{C}M r_{k+1}^{\beta}\left(\frac{1-\mu}{\tau_1^{\beta }}+\frac{\mu}{\hat{C}\tau_1^{\beta }}+\frac{C}{\hat{C}\tau_1^{\beta }}\right).
\end{aligned}
\end{equation*}

Take $\beta$ small enough and $\hat{C}$ large enough such that
\begin{equation*}
\frac{1-\mu}{\tau_1^{\beta }}+\frac{\mu}{\hat{C}\tau_1^{\beta }}+\frac{C}{\hat{C}\tau_1^{\beta }}\leq 1.
\end{equation*}
Then
\begin{equation*}
\sup_{\Omega _{r_{k+1}}} u\leq \sup_{\Omega _{r_{k+1}}} v\leq \hat{C} M r_{k+1}^{\beta}.
\end{equation*}
The proof for
\begin{equation*}
\inf_{\Omega _{r_{k+1}}} u\geq -\hat{C} M r_{k+1}^{\beta}
\end{equation*}
is similar and we omit here. Therefore,
\begin{equation*}
\|u\|_{L^{\infty}(\Omega _{r_{k+1}})}\leq \hat{C} M r_{k+1}^{\beta}.
\end{equation*}
Hence,\cref{e-poisson-discrete} holds for $k+1$. By induction, the proof is completed. \qed~\\

\begin{remark}\label{r-1}
In above theorem, we do not indicate which kind of solutions we consider. In fact, the key points are the strong maximum principle in proving an estimate similar to\cref{e.2} (see the following proofs for details), and the scaling. Hence, the solution can be understood in any sense (such as in the sense of weak solution or strong solution etc.) only if the strong maximum principle holds.
\end{remark}

In the following, we will prove the boundary H\"{o}lder regularity for general elliptic equations under the (H2) condition. First, we introduce a sequence functions based on the (H2) condition. Choose and fix one $g_{\nu}\in C^{\infty}(\partial B_1)$ with $0\leq g_{\nu} \leq 1$ and
\begin{equation}\label{gnu}
g_{\nu}(x)\equiv
\left\{ \begin{aligned}
&0 ~~&&\mbox{on}~~\partial B_1\cap B(e_1,\nu/2);\\
&1 ~~&&\mbox{on}~~\partial B_1\backslash B(e_1,\nu),
\end{aligned}\right.
\end{equation}
where $e_1=(1,0,0,...,0)$. Clearly, $g_{\nu}$ depends only on $n$ and $\nu$. Next, for each $k\geq 1$, define
\begin{equation*}\label{gk3}
\tilde{g}_k(x)=g_{\nu}(x/r_k),
\end{equation*}
which is defined on $\partial B_{r_k}$. Since $y_k\in \partial B_{r_k}$, there exists an orthogonal matrix $T_k$ such that $r_ke_1=T_k y_k$. Finally, introduce
\begin{equation}\label{gk2}
g_k(x)= \tilde{g}_{k}(T_kx)=g_{\nu}(T_kx/r_k).
\end{equation}
Then
\begin{equation*}\label{gk1}
g_{k}(x)\equiv
\left\{ \begin{aligned}
&0 ~~&&\mbox{on}~~\partial B_{r_k}\cap B(y_k,\nu r_k/2);\\
&1 ~~&&\mbox{on}~~\partial B_{r_k}\backslash B(y_k,\nu r_k).
\end{aligned}\right.
\end{equation*}
In above definition, $\nu,\{r_k\}$ and $\{y_k\}$ are as in the (H2) condition.

First, we consider the uniformly elliptic equations in divergence form. We begin with the following lemma.
\begin{lemma}\label{l-div}
Let $a^{ij}$ be uniformly elliptic with constants $\lambda$ and $\Lambda$, and $v$ be a weak solution of
\begin{equation*}
\left\{\begin{aligned}
&(a^{ij} v_{i})_{j}=0 ~~&&\mbox{in}~~B_1;\\
& v=g ~~&&\mbox{on}~~\partial B_1,
\end{aligned}
\right.
\end{equation*}
where $g\in C^{\infty}(\partial B_1)$, $0\leq g \leq 1$ and $g\equiv 0$ on a portion of $\partial B_1$. Then for any $0<\delta<1$,
\begin{equation*}\label{e-div-l-decay}
  \sup_{B_{\delta}}v\leq 1-\mu,
\end{equation*}
where $0<\mu<1$ depends only on $n,\lambda,\Lambda,\delta$ and $g$.
\end{lemma}
\proof By the global H\"{o}lder estimate (see \cite[Theorem 8.29]{MR1814364}),
\begin{equation*}
  \|v\|_{C^{\gamma}(\bar{B}_{1})}\leq C_0 \|g\|_{C^{1}(\bar{B}_{1})},
\end{equation*}
where $0<\gamma<1$ and $C_0$ depend only on $n,\lambda$ and $\Lambda$. Let $x_0\in \partial B_1$ with $g(x_0)=0$. Take $0<t<1-\delta$ small enough (depending only on $n,\lambda,\Lambda,\delta$ and $g$) such that
\begin{equation*}
  v((1-t)x_0)=v((1-t)x_0)-v(x_0)\leq t^{\gamma}C_0\|g\|_{C^{1}(\bar{B}_{1})}\leq 1/2.
\end{equation*}
By the interior Harnack inequality (see \cite[Theorem 8.20]{MR1814364}),
\begin{equation*}
  1- v(x)\geq c_0\left(1- v((1-t)x_0)\right)\geq c_0/2, ~~\forall ~x\in \partial B_{1-t},
\end{equation*}
where $c_0$ depends only on $n,\lambda,\Lambda$ and $t$. Hence,
\begin{equation*}
  \sup_{B_{\delta}} v\leq \sup_{B_{1-t}} v \leq 1-c_0/2:=1-\mu.
\end{equation*}
\qed~\\

Based on above result, the boundary H\"{o}lder regularity can be derived for linear elliptic equations in divergence form:
\begin{equation}\label{e-div-2}
 \left\{ \begin{aligned}
&(a^{ij}u_i)_j=f ~~&&\mbox{in}~~\Omega;\\
&u=g ~~&&\mbox{on}~~\partial \Omega.
 \end{aligned}\right.
\end{equation}

\begin{theorem}\label{t-div-2}
Suppose that $\Omega$ satisfies the (H2) condition at $0\in \partial \Omega$. Let $u$ be a weak solution of\cref{e-div-2} where $a^{ij}$ is uniformly elliptic with $\lambda$ and $\Lambda$, $f\in L^{p}(\Omega)$ with $p>n/2$ and $g\in C^{\alpha}(0)$. Then $u$ is $C^{\beta}$ at $0$ and
\begin{equation*}
  |u(x)-u(0)|\leq C |x|^{\beta}\left(\|u\|_{L^{\infty }(\Omega\cap B_1)}+\|f\|_{L^{p}(\Omega\cap B_1)}+[g]_{C^{\alpha}(0)}\right), ~~\forall ~x\in \Omega\cap B_1,
\end{equation*}
where $0<\beta\leq \min(2-n/p,\alpha)$ depends only on $n,\lambda,\Lambda, \tau_1, \tau_2,\nu$ and $p$; $C$ depends only on $n,\lambda,\Lambda,\tau_2,\nu$ and $p$. Here, $\tau_1, \tau_2$ and $\nu$ are constants from \Cref{d-1.1} and the (H2) condition.
\end{theorem}
\proof We assume that $g(0)=0$ as before. Let $M=\|u\|_{L^{\infty }(\Omega\cap B_1)}+\|f\|_{L^{p}(\Omega\cap B_1)}+[g]_{C^{\alpha}(0)}$ and $\Omega_{r}=\Omega\cap B_{r}$. We only need to prove the following:

There exist constants $0<\beta\leq \min(2-n/p,\alpha)$ depending only on $n,\lambda,\Lambda, \tau_1, \tau_2,\nu$ and $p$, and $\hat{C}$ depending only on $n,\lambda,\Lambda,\tau_2,\nu$ and $p$ such that for all $k\geq 0$,

\begin{equation}\label{e-div-discrete}
\|u\|_{L^{\infty }(\Omega _{r_{k}})}\leq \hat{C} M r_{k}^{\beta}.
\end{equation}

We prove\cref{e-div-discrete} by induction. For $k=0$, it holds clearly. Suppose that it holds for $k$. We need to prove that it holds for $k+1$.

Let $v$ solve (see \cite[Theorem 5.2]{MR0244627})
\begin{equation*}
\left\{\begin{aligned}
 &(a^{ij}v_i)_j=0 &&\mbox{in}~~B_{r_k}; \\
 &v=\tilde{g} &&\mbox{on}~~\partial B_{r_k},
\end{aligned}
\right.
\end{equation*}
where $\tilde{g}=(\hat{C} M r_{k}^{\beta}-Mr_{k}^{\alpha})g_k+Mr_{k}^{\alpha}$ and $g_k$ is defined in\cref{gk2}. Let $w$ solve (see also \cite[Theorem 5.2]{MR0244627})
\begin{equation*}
\left\{\begin{aligned}
 &(a^{ij}w_i)_j=-|f| &&\mbox{in}~~B_{r_k}; \\
 &w=0 &&\mbox{on}~~\partial B_{r_k}.
\end{aligned}
\right.
\end{equation*}
Then
\begin{equation*}
 (a^{ij}(u-v-w)_i)_j\geq 0~~ \mbox{in}~~\Omega_{r_k}.
\end{equation*}

Note that $g_k\geq 0$ on $\partial B_{r_k}$ and $w\geq 0$ in $B_{r_k}$. Thus,
\begin{equation*}
  u=g\leq Mr_k^{\alpha}=\min_{\partial B_{r_k}} \tilde{g}\leq v\leq v+w~~\mbox{on}~~\partial \Omega\cap \bar{B}_{r_k}.
\end{equation*}
On the other hand, since $g_k\geq 1$ on $\partial B_{r_k}\cap \Omega$,
\begin{equation*}
  u\leq \hat{C}Mr_k^{\beta}\leq \tilde{g}~~\mbox{on}~~\partial B_{r_k}\cap \Omega.
\end{equation*}
Hence,
\begin{equation*}
  u\leq v+w~~\mbox{on}~~\partial \Omega_{r_k}.
\end{equation*}
Similarly, we have
\begin{equation*}
  u\geq -(v+w)~~\mbox{on}~~\partial \Omega_{r_k}.
\end{equation*}
Then by the maximum principle (see \cite[Theorem 8.1]{MR1814364}),
\begin{equation}\label{e-div-6}
  -v-w\leq u \leq v+w ~~\mbox{in}~~\Omega_{r_k}.
\end{equation}

Now, we estimate $v$ and $w$ respectively. Let $y=T_kx/r_k$ and
\begin{equation*}
  \tilde{v}(y)=\frac{v(x)- Mr_k^{\alpha}}{\hat{C}Mr_k^{\beta}-Mr_k^{\alpha}}.
\end{equation*}
Then $\tilde{v}$ satisfy
\begin{equation*}
\left\{\begin{aligned}
 &(\tilde{a}^{ij}\tilde{v}_i)_j=0 &&\mbox{in}~~B_{1}; \\
 &\tilde{v}=g_{\nu} &&\mbox{on}~~\partial B_{1},
\end{aligned}
\right.
\end{equation*}
where $\tilde{a}^{ij}(y)=a^{ij}(x)T_k^{li}T_k^{mj}$ is uniformly elliptic with $\lambda$ and $\Lambda$. By \Cref{l-div} we have
\begin{equation*}\label{e-div-4-2}
  \sup_{B_{\tau_2}}v\leq 1-\mu,
\end{equation*}
where $0<\mu<1$ depends only on $n,\lambda,\Lambda,\tau_2$ and $g_{\nu}$. Since $g_{\nu}$ is determined by $\nu$ (see\cref{gnu}),  $\mu$ depends only on $n,\lambda,\Lambda,\tau_2$ and $\nu$. Rescaling back to $v$, we have
\begin{equation}\label{e-div-4}
\begin{aligned}
\sup_{B_{r_{k+1}}}v\leq\sup_{B_{\tau_2r_{k}}}v&\leq (1-\mu)\left(\hat{C}M r_{k}^{\beta}-M r_{k}^{\alpha}\right)+M r_{k}^{\alpha}\\
&\leq \hat{C}M r_{k+1}^{\beta}\cdot \frac{1-\mu}{\tau_1^{\beta }}+\mu M r_{k}^{\beta}\\
&\leq \hat{C}M r_{k+1}^{\beta}\left(\frac{1-\mu}{\tau_1^{\beta }}+\frac{\mu}{\hat{C}\tau_1^{\beta }}\right).
\end{aligned}
\end{equation}

By the Alexandrov-Bakel'man-Pucci maximum principle for $w$ (see \cite[Theorem 8.16]{MR1814364}), we have
\begin{equation}\label{e-div-5}
\begin{aligned}
  \|w\|_{L^{\infty}(B_{r_{k+1}})}&\leq CMr_{k}^{2-n/p}\leq \hat{C}M r_{k+1}^{\beta}\cdot\frac{C}{\hat{C}\tau_1^{\beta }},
\end{aligned}
\end{equation}
where $C$ depends only on $n, \lambda, \Lambda$ and $p$.

Take $\beta$ small enough and $\hat{C}$ large enough such that
\begin{equation*}
\frac{1-\mu}{\tau_1^{\beta }}+\frac{\mu}{\hat{C}\tau_1^{\beta }}+\frac{C}{\hat{C}\tau_1^{\beta }}\leq 1.
\end{equation*}
Then combining\cref{e-div-6},\cref{e-div-4} and\cref{e-div-5}, we have
\begin{equation*}
\|u\|_{L^{\infty }(\Omega _{r_{k+1}})}\leq \hat{C} M r_{k+1}^{\beta}.
\end{equation*}
By induction, the proof is completed. \qed~\\

\begin{remark}
In above proof, the only used tools are the solvability of $v$ in a ball, the strong maximum principle for $v$, the comparison principle between $u$ and $v$ and a scaling argument. These are valid for many types of equations. The benefit of the method is that it does not need to construct directly a barrier which is difficult for domains with complicated geometric structures.
\end{remark}

Our method is also applicable for elliptic equations in non-divergence form. More generally, we consider the fully nonlinear elliptic equations in non-divergence form:
\begin{equation}\label{e-ful-1}
\left\{\begin{aligned}
  &u\in S(\lambda,\Lambda,f) &&~~\mbox{in}~~\Omega;\\
  &u=g &&~~\mbox{on}~~\partial \Omega.
\end{aligned}\right.
\end{equation}
Here, we consider the viscosity solutions and the notations are adopted from \cite{MR1351007} (see also \cite{MR1376656} and \cite{MR1118699}).
\begin{theorem}\label{t-full-elliptic}
Suppose that $\Omega$ satisfies the (H2) condition at $0\in \partial \Omega$. Let $u$ be a viscosity solution of\cref{e-ful-1} where $f\in L^{n}(\Omega)$ and $g\in C^{\alpha}(0)$. Then $u$ is $C^{\beta}$ at $0$ and
\begin{equation*}
  |u(x)-u(0)|\leq C |x|^{\beta}\left(\|u\|_{L^{\infty }(\Omega\cap B_1)}+\|f\|_{L^{n}(\Omega\cap B_1)}+[g]_{C^{\alpha}(0)}\right), ~~\forall ~x\in \Omega\cap B_1,
\end{equation*}
where $0<\beta\leq \alpha$ depends only on $n, \lambda, \Lambda, \tau_1, \tau_2$ and $\nu$; $C$ depend only on $n, \lambda, \Lambda, \tau_2$ and $\nu$. Here, $\tau_1, \tau_2$ and $\nu$ are constants from \Cref{d-1.1} and the (H2) condition.
\end{theorem}
\proof We assume that $g(0)=0$ as before. Let $M=\|u\|_{L^{\infty }(\Omega\cap B_1)}+\|f\|_{L^{n}(\Omega\cap B_1)}+[g]_{C^{\alpha}(0)}$ and $\Omega_{r}=\Omega\cap B_{r}$. We only need to prove the following:

There exist constants $0<\beta\leq \alpha$ depending only on $n, \lambda, \Lambda, \tau_1, \tau_2$ and $\nu$, and $\hat{C}$ depending only on $n,\lambda,\Lambda,\tau_2$ and $\nu$ such that for all $k\geq 0$,

\begin{equation}\label{e-ful-discrete}
\|u\|_{L^{\infty }(\Omega _{r_{k}})}\leq \hat{C} M r_{k}^{\beta}.
\end{equation}

We prove\cref{e-ful-discrete} by induction. For $k=0$, it holds clearly. Suppose that it holds for $k$. We need to prove that it holds for $k+1$.

Let $v$ solve (see \cite[Corollary 3.10]{MR1376656})
\begin{equation*}
\left\{\begin{aligned}
 &M^{+}(D^2v,\lambda,\Lambda)=0 &&\mbox{in}~~B_{r_k}; \\
 &v=\tilde{g} &&\mbox{on}~~\partial B_{r_k},
\end{aligned}
\right.
\end{equation*}
where $\tilde{g}=(\hat{C} M r_{k}^{\beta}-Mr_{k}^{\alpha})g_k+Mr_{k}^{\alpha}$. Let $w$ solve (see also \cite[Corollary 3.10]{MR1376656})
\begin{equation*}
\left\{\begin{aligned}
 &M^{+}(D^2w,\lambda,\Lambda)=-|f| &&\mbox{in}~~B_{r_k}; \\
 &w=0 &&\mbox{on}~~\partial B_{r_k}.
\end{aligned}
\right.
\end{equation*}
Then
\begin{equation*}
M^{+}(D^2(v+w),\lambda,\Lambda)\leq f ~~\mbox{in}~~\Omega_{r_k}.
\end{equation*}
As in the proof of \Cref{t-div-2}, it can be verified similarly that
\begin{equation*}
u\leq v+w~\mbox{on}~~\partial \Omega_{r_k}.
\end{equation*}
Hence, by the comparison principle (see \cite[Theorem 3.3]{MR1118699}), we have
\begin{equation}\label{e-ful-2}
u \leq v+w ~~\mbox{in}~~\Omega_{r_k}.
\end{equation}

Similar to the proof of \Cref{t-div-2}, we estimate $v$ and $w$ respectively. Let $y=T_kx/r_k$ and
\begin{equation*}
  \tilde{v}(y)=\frac{v(x)- Mr_k^{\alpha}}{\hat{C}Mr_k^{\beta}-Mr_k^{\alpha}}.
\end{equation*}
Then $\tilde{v}$ satisfy
\begin{equation*}
\left\{\begin{aligned}
 &M^{+}(D^2\tilde{v},\lambda,\Lambda)=0 &&\mbox{in}~~B_{1}; \\
 &\tilde{v}=g_{\nu} &&\mbox{on}~~\partial B_{1}.
\end{aligned}
\right.
\end{equation*}
By the strong maximum principle (see \cite[Proposition 4.9]{MR1351007}), we have
\begin{equation*}
  \sup_{B_{\tau_2}} \tilde{v} \leq 1-\mu,
\end{equation*}
where $0<\mu<1$. Note that $\mu$ is determined obviously by $\tau_2$ and $\tilde{v}$, and the later is uniquely determined by the operator $M^+$, the domain $B_1$ and the boundary value $g_{\nu}$. Since $M^+$ depends only on $\lambda$ and $\Lambda$, $B_1$ depends only on $n$ and $g_{\nu}$ depends only on $n$ and $\nu$, $\mu$ depends only on $n,\lambda,\Lambda,\tau_2$ and $\nu$.

Rescaling back to $v$, we have
\begin{equation}\label{e-ful-3}
\begin{aligned}
\sup_{B_{r_{k+1}}}v\leq\sup_{B_{\tau_2r_{k}}}v&\leq (1-\mu)\left(\hat{C}M r_{k}^{\beta}-M r_{k}^{\alpha}\right)+M r_{k}^{\alpha}\\
&\leq \hat{C}M r_{k+1}^{\beta}\cdot \frac{1-\mu}{\tau_1^{\beta }}+\mu M r_{k}^{\beta}\\
&\leq \hat{C}M r_{k+1}^{\beta}\left(\frac{1-\mu}{\tau_1^{\beta }}+\frac{\mu}{\hat{C}\tau_1^{\beta }}\right).
\end{aligned}
\end{equation}

By the Alexandrov-Bakel'man-Pucci maximum principle for $w$ (see \cite[Theorem 3.2]{MR1351007}), we have
\begin{equation}\label{e-ful-4}
\begin{aligned}
  \|w\|_{L^{\infty}(B_{r_{k+1}})}\leq CMr_{k} \leq \hat{C}M r_{k+1}^{\beta}\cdot\frac{C}{\hat{C}\tau_1^{\beta }},
\end{aligned}
\end{equation}
where $C$ depends only on $n, \lambda$ and $\Lambda$.

Take $\beta$ small enough and $\hat{C}$ large enough such that
\begin{equation*}
\frac{1-\mu}{\tau_1^{\beta }}+\frac{\mu}{\hat{C}\tau_1^{\beta }}+\frac{C}{\hat{C}\tau_1^{\beta }}\leq 1.
\end{equation*}
Then combining\cref{e-ful-2},\cref{e-ful-3} and\cref{e-ful-4}, we have
\begin{equation*}
\sup_{\Omega _{r_{k+1}}} u\leq \hat{C} M r_{k+1}^{\beta}.
\end{equation*}
The proof for
\begin{equation*}
\inf_{\Omega _{r_{k+1}}} u\geq -\hat{C} M r_{k+1}^{\beta}
\end{equation*}
is similar and we omit here. Therefore,
\begin{equation*}
\|u\|_{L^{\infty}(\Omega _{r_{k+1}})}\leq \hat{C} M r_{k+1}^{\beta}.
\end{equation*}
By induction, the proof is completed. \qed~\\

Quasilinear elliptic equations can also be treated. We prove the boundary H\"{o}lder regularity for $p-$Laplace equations:
\begin{equation}\label{e-deg-1}
\left\{\begin{aligned}
 \mathrm{div} (|\nabla u|^{p-2} \nabla u)&=0&&~~\mbox{in}~~\Omega;\\
  u&=g &&~~\mbox{on}~~\partial \Omega.
\end{aligned}\right.
\end{equation}

\begin{theorem}\label{t-deg-elliptic}
Suppose that $\Omega$ satisfies the (H2) condition at $0\in \partial \Omega$. Let $u$ be a weak solution of\cref{e-deg-1} where $1< p<+\infty$ and $g\in C^{\alpha}(0)$. Then $u$ is $C^{\beta}$ at $0$ and
\begin{equation*}
  |u(x)-u(0)|\leq 8|x|^{\beta}\left(\|u\|_{L^{\infty }(\Omega\cap B_1)}+[g]_{C^{\alpha}(0)}\right), ~~\forall ~x\in \Omega\cap B_1,
\end{equation*}
where $0<\beta\leq \alpha$ depends only on $n,\tau_1,\tau_2,\nu$ and $p$. Here, $\tau_1, \tau_2$ and $\nu$ are constants from \Cref{d-1.1} and the (H2) condition.
\end{theorem}
\proof We assume that $g(0)=0$ as before. Let $M=\|u\|_{L^{\infty }(\Omega\cap B_1)}+[g]_{C^{\alpha}(0)}$ and $\Omega_{r}=\Omega\cap B_{r}$. We only need to prove the following:

There exist a constant $0<\beta\leq \alpha$ depending only on $n,\tau_1,\tau_2,\nu$ and $p$ such that
\begin{equation}\label{e-p-place-1}
  \tau_1^{\beta}\geq \frac{1}{2}
\end{equation}
and for all $k\geq 0$,
\begin{equation}\label{e-deg-discrete}
\|u\|_{L^{\infty }(\Omega _{r_{k}})}\leq 4M r_{k}^{\beta}.
\end{equation}

We prove\cref{e-deg-discrete} by induction. For $k=0$, it holds clearly. Suppose that it holds for $k$. We need to prove that it holds for $k+1$.

Let $v$ solve (see \cite[Theorem 2.16]{Lindqvist})
\begin{equation*}\label{e-deg-2}
\left\{\begin{aligned}
 \mathrm{div} (|\nabla v|^{p-2} \nabla v)&=0 &&\mbox{in}~~B_{r_k}; \\
 v&=\tilde{g} &&\mbox{on}~~\partial B_{r_k},
\end{aligned}
\right.
\end{equation*}
where $\tilde{g}=(4 M r_{k}^{\beta}-Mr_{k}^{\alpha})g_k+Mr_{k}^{\alpha}$. As in the proof of \Cref{t-div-2}, it can be verified similarly that
\begin{equation*}
-v\leq u\leq v~\mbox{on}~~\partial \Omega_{r_k}.
\end{equation*}
Hence, by the comparison principle (see \cite[Theorem 2.16]{Lindqvist}), we have
\begin{equation}\label{e-deg-3}
-v\leq u \leq v ~~\mbox{in}~~\Omega_{r_k}.
\end{equation}
Similarly to the previous proof, let $y=T_kx/r_k$ and
\begin{equation*}
  \tilde{v}(y)=\frac{v(x)- Mr_k^{\alpha}}{\hat{C}Mr_k^{\beta}-Mr_k^{\alpha}}.
\end{equation*}
Then $\tilde{v}$ satisfy
\begin{equation*}
\left\{\begin{aligned}
 \mathrm{div} (|\nabla \tilde{v}|^{p-2} \nabla \tilde{v})&=0 &&\mbox{in}~~B_{1}; \\
 \tilde{v}&=g_{\nu} &&\mbox{on}~~\partial B_{1}.
\end{aligned}
\right.
\end{equation*}

By the strong maximum principle (see \cite[Corollary 2.21]{Lindqvist}), we have
\begin{equation*}
  \sup_{B_{\tau_2}} \tilde{v} \leq 1-\mu,
\end{equation*}
where $0<\mu<1$. Note that $\mu$ is determined obviously by $\tau_2$ and $\tilde{v}$, and the later is uniquely determined by the $p-$Laplace operator (depending only on $p$), the domain $B_1$ (depending only on $n$) and the boundary value $g_{\nu}$ (depending only on $n$ and $\nu$). Hence, $\mu$ depends only on $n,\tau_2,\nu$ and $p$.

Rescaling back to $v$, we have
\begin{equation}\label{e-deg-4}
\begin{aligned}
\sup_{B_{r_{k+1}}}v\leq\sup_{B_{\tau_2r_{k}}}v&\leq (1-\mu)\left(4M r_{k}^{\beta}-M r_{k}^{\alpha}\right)+M r_{k}^{\alpha}\\
&\leq 4M r_{k+1}^{\beta}\cdot \frac{1-\mu}{\tau_1^{\beta }}+\mu M r_{k}^{\beta}\\
&\leq 4M r_{k+1}^{\beta}\left(\frac{1-\mu}{\tau_1^{\beta }}+\frac{\mu}{4\tau_1^{\beta }}\right).
\end{aligned}
\end{equation}

Take $\beta$ small enough such that\cref{e-p-place-1} holds and
\begin{equation*}
\frac{1-\mu}{\tau_1^{\beta }}<1-\frac{\mu}{2}.
\end{equation*}
Then combining\cref{e-deg-3} and\cref{e-deg-4}, we have
\begin{equation*}
\|u\|_{L^{\infty}(\Omega _{r_{k+1}})}\leq 4 M r_{k+1}^{\beta}.
\end{equation*}
By induction, the proof is completed. \qed~\\

\begin{remark}\label{r-2}
Although we apply weak solutions in \Cref{t-div-2} and \Cref{t-deg-elliptic}, and viscosity solutions in \Cref{t-full-elliptic}, the kinds of solutions do not play a key role in the theory since the critical elements are the strong maximum principle and the scaling.
\end{remark}

Finally, we consider the fractional Laplace equations:
\begin{equation}\label{e-frac-1}
\left\{\begin{aligned}
 (-\Delta )^{s/2}u&=f ~~\mbox{in}~~\Omega; \\
 u&=g~~\mbox{in}~~\Omega^c,
\end{aligned}
\right.
\end{equation}
where $0<s<2$. For the fundamental theory of the fractional Laplace equations, we refer to \cite{C-L-M}. We have the following boundary H\"{o}lder regularity corresponding to\cref{e-frac-1}.
\begin{theorem}\label{t-frac}
Suppose that $\Omega$ satisfies the (H3) condition at $0\in \partial \Omega$. Let $u$ satisfy\cref{e-frac-1} where $f\in L^{p}(\Omega)$ with $p>n/s$, and $g$ is bounded and $C^{\alpha}$ at $0$. Then $u$ is $C^{\beta}$ at $0$ and
\begin{equation*}
  |u(x)-u(0)|\leq C |x|^{\beta}\cdot\left(\|u\|_{L^{\infty}(\Omega\cap B_1)}+\|f\|_{L^{p}(\Omega\cap B_1)}+
  \|g\|_{C^{\alpha}(0)}\right), ~~\forall ~x\in \Omega\cap B_1,
\end{equation*}
where $0<\beta\leq \min(s-n/p,\alpha)$ depends only on $n,s,\tau_1, \tau_2$ and $\nu$; $C$ depends only on $n,s,\tau_1,\tau_2,\nu$ and $p$. Here, $\tau_1, \tau_2$ and $\nu$ are constants from \Cref{d-1.1} and the (H3) condition.
\end{theorem}

\begin{remark}
Ros-Oton and Serra \cite{MR3168912} proved the boundary H\"{o}lder regularity for\cref{e-frac-1} with $f\in L^{\infty}$, $g\equiv 0$ and the exterior sphere condition at $0$.
\end{remark}

\noindent\emph{Proof of \Cref{t-frac}. }The essential of the proof is the same as previous and we also assume that $g(0)=0$. Let $M=\|u\|_{L^{\infty}(\Omega\cap B_1)}+\|f\|_{L^{p}(\Omega\cap B_1)}+\|g\|_{C^{\alpha}(0)}$. We only need to prove the following:

There exist constants $0< \beta\leq \min(s-n/p,\alpha)$ and $k_0\geq 1$ depending only on $n,s,\tau_1, \tau_2$ and $\nu$, $\hat{C}$ depending only on $n,s,\tau_1,\tau_2,\nu$ and $p$, and a sequence of nonnegative functions $\{v_k\}_{k=1}^{\infty}$ with $v_1\equiv M$ such that for $k\geq 2$,
\begin{equation}\label{e-final-3}
\left\{
\begin{aligned}
&(-\Delta)^{s/2} v_k=|f| &&~~\mbox{in}~~B_{r_{(k-1)k_0}}; \\
&-v_k \leq u \leq v_k &&~~\mbox{in}~~R^{n}; \\
&v_k \leq \hat{C}Mr_{kk_0}^{\beta}&&~~\mbox{in}~~B_{r_{kk_0}}.
\end{aligned}\right.
\end{equation}

Indeed, by\cref{e-final-3}, we have
\begin{equation*}\label{e-frac-discrete}
\|u\|_{L^{\infty }(B_{r_{kk_0}})}\leq \hat{C} M r_{kk_0}^{\beta},
\end{equation*}
which implies the boundary H\"{o}lder regularity in the same way as previous.

We prove\cref{e-final-3} by induction. For $k=2$, let
\begin{equation*}
v_{2}(x)=\int_{B^c_{r_{k_0}}} P_{r_{k_0}}(x,y)v_1(y) dy+\int_{ B_{r_{k_0}}} G(x,y)|f|(y)  dy,
\end{equation*}
where $x\in B_{r_{k_0}}$, $P_r$ is the Poisson kernel (see \cite[Chap. 4.1]{C-L-M}) for $B_r$ with
\begin{equation*}
P_r(x,y)=\frac{\Gamma(n/2)}{\pi^{n/2+1}}\cdot\sin \frac{\pi s}{2} \cdot \left(\frac{r^2-|x|^2}{|y|^2-r^2}\right)^{s/2}\cdot\frac{1}{|x-y|^n},~~x\in B_r, y\in B_r^c,
\end{equation*}
and $G(x,y)$ is the Green function (see \cite[Chap. 2.2]{C-L-M}) with
\begin{equation*}
  G(x,y)=\frac{C(n,s)}{|x-y|^{n-s}}+\frac{c(n,s)}{|x-y|^{n-s}}
\end{equation*}
for some positive constants $C(n,s)$ and $c(n,s)$ depending only on $n$ and $s$. Then $v$ satisfies (see \cite[Chap. 4.1]{C-L-M})
\begin{equation*}
\left\{
\begin{aligned}
&(-\Delta)^{s/2} v_2=|f| &&~~\mbox{in}~~B_{r_{k_0}}; \\
&v_2=M&&~~\mbox{in}~~B^c_{r_{k_0}}.
\end{aligned}\right.
\end{equation*}

In aid of the maximum principle, the first two terms in\cref{e-final-3} hold clearly. Since $P_{r_{k_0}}$ is the Poisson kernel,
\begin{equation*}
\int_{B^c_{r_{k_0}}} P_{r_{k_0}}(x,y)v_1(y) dy=\int_{B^c_{r_{k_0}}} P_{r_{k_0}}(x,y)M dy\equiv M.
\end{equation*}
By H\"{o}lder's inequality,
\begin{equation}\label{e-final-4}
\begin{aligned}
\int_{ B_{r_{k_0}}} G(x,y)|f|(y)  dy\leq C_1Mr_{k_0}^{s-n/p}
\end{aligned}
\end{equation}
for some constant $C_1$ depending only on $n, s$ and $p$. Hence,
\begin{equation}\label{e-final-5}
\begin{aligned}
  v_{2}(x)&\leq M+C_1Mr_{k_0}^{s-n/p} \\
  &\leq M+C_1Mr_{k_0}^{\beta}\\
  &\leq \hat{C}Mr_{2k_0}^{\beta}\left(\frac{1}{\hat{C}\tau_1^{2k_0\beta}}
  +\frac{C_1}{\hat{C}\tau_1^{k_0\beta}}\right).
\end{aligned}
\end{equation}

Take proper $k_0$ and $\beta$ (to be specified later) such that
\begin{equation}\label{e-final-6}
  \tau_1^{2k_0\beta}\geq 1/2.
\end{equation}
Next, take $\hat{C}$ large enough such that
\begin{equation}\label{e-final-10}
\frac{1}{\hat{C}\tau_1^{2k_0\beta}}
  +\frac{C_1}{\hat{C}\tau_1^{k_0\beta}}\leq 1.
\end{equation}
Then we have
\begin{equation*}\label{e-Lap-decay-1}
\begin{aligned}
\sup_{B_{r_{2k_0}}} v_2 &\leq \hat{C}M r_{2k_0}^{\beta}.
\end{aligned}
\end{equation*}
That is, the third term in\cref{e-final-3} holds. Hence,\cref{e-final-3} holds for $k=2$.

Suppose that\cref{e-final-3} holds for $k$. We need to prove that\cref{e-final-3} holds for $k+1$. Let
\begin{equation*}
  \tilde{g}(x)\equiv  \left\{
\begin{aligned}
 &M r_{(k-1)k_0}^{\alpha},&&~~\mbox{in}~~ (B_{r_{(k-1)k_0}}\backslash B_{r_{kk_0}})\cap \Omega^c; \\
 &\hat{C}M r_{(k-1)k_0}^{\beta},&&~~\mbox{in}~~ (B_{r_{(k-1)k_0}}\backslash B_{r_{kk_0}})\cap \Omega;\\
 &v_k,&&~~\mbox{in}~~ R^{n}\backslash B_{r_{(k-1)k_0}}.
\end{aligned}\right.
\end{equation*}
Then, it is easy to check that $-\tilde{g}\leq u\leq \tilde{g}$ in $B^c_{r_{kk_0}}$. Let
\begin{equation}\label{e-final-v}
v_{k+1}(x)=\int_{B^c_{r_{kk_0}}} P_{r_{kk_0}}(x,y)\tilde{g}(y) dy +\int_{ B_{r_{kk_0}}} G(x,y)|f|(y)  dy,
\end{equation}
where $x\in B_{r_{kk_0}}$. Then the first two terms in\cref{e-final-3} hold clearly.

Split the first term in\cref{e-final-v} into two parts:
\begin{equation}\label{e-final-11}
\begin{aligned}
&\int_{B^c_{r_{kk_0}}} P_{r_{kk_0}}\tilde{g}dy=\int_{B_{r_{(k-1)k_0}}\backslash B_{r_{kk_0}}} P_{r_{kk_0}}\tilde{g}dy+\int_{B^c_{r_{(k-1)k_0}}} P_{r_{kk_0}}\tilde{g} dy.
\end{aligned}
\end{equation}
Now, we estimate them separately. Let $A=(B_{r_{kk_0-1}}\backslash B_{r_{kk_0}})\cap \Omega^c,\tilde{x}=x/r_{kk_0},\tilde{y}=y/r_{kk_0}$ and $\tilde{A}=A/r_{kk_0}$. For $x\in B_{r_{kk_0}}$,
\begin{equation*}\label{e-final-7-2}
\begin{aligned}
&\int_{A} P_{r_{kk_0}}(x,y) dy\\
&=C\int_{A} \left(\frac{r_{kk_0}^2-|x|^2}{|y|^2-r_{kk_0}^2}\right)^{s/2}\cdot\frac{1}{|x-y|^n} dy\\
&=C\int_{\tilde{A}} \left(\frac{1-|\tilde{x}|^2}{|\tilde{y}|^2-1}\right)^{s/2}\cdot\frac{1}{|\tilde{x}-\tilde{y}|^n} d\tilde{y},
\end{aligned}
\end{equation*}
where $C$ depends only on $n$ and $s$. By the (H3) condition, $|\tilde{A}|\geq \nu$. Since $\tilde{x}\in B_1$, $\tilde{y}\in B_1/\tau_1$ and the integrand is positive,
\begin{equation}\label{e-final-7-3}
\int_{A} P_{r_{kk_0}}(x,y) dy=C\int_{\tilde{A}} \left(\frac{1-|\tilde{x}|^2}{|\tilde{y}|^2-1}\right)^{s/2}\cdot\frac{1}{|\tilde{x}-\tilde{y}|^n} d\tilde{y}\geq \mu,
\end{equation}
where $\mu$ depends only on $n,s,\tau_1,\tau_2$ and $\nu$. Let
\begin{equation*}
  \hat{g}=\frac{\tilde{g}-M r_{(k-1)k_0}^{\alpha}}{\hat{C}M r_{(k-1)k_0}^{\beta}-M r_{(k-1)k_0}^{\alpha}}.
\end{equation*}
Then $ \hat{g}\leq 1$ in $B_{r_{(k-1)k_0}}$ and $\hat{g}\leq 0$ in $B_{r_{(k-1)k_0}}\cap \Omega^c$. Thus,
\begin{equation*}
\begin{aligned}
&\int_{B_{r_{(k-1)k_0}}\backslash B_{r_{kk_0}}} P_{r_{kk_0}}(x,y)\hat{g}(y) dy\\
&=\int_{B_{r_{(k-1)k_0}}\backslash B_{r_{kk_0}}} P_{r_{kk_0}}(x,y) dy+\int_{B_{r_{(k-1)k_0}}\backslash B_{r_{kk_0}}} P_{r_{kk_0}}(x,y)(\hat{g}(y)-1) dy\\
&\leq \int_{B^c_{r_{kk_0}}} P_{r_{kk_0}}(x,y) dy-\int_{(B_{r_{(k-1)k_0}}\backslash B_{r_{kk_0}})\cap\Omega^c} P_{r_{kk_0}}(x,y) dy\\
&\leq 1-\int_{(B_{r_{kk_0-1}}\backslash B_{r_{kk_0}})\cap\Omega^c} P_{r_{kk_0}}(x,y) dy\\
&\leq 1-\mu ~~\mbox{by}\cref{e-final-7-3}.
\end{aligned}
\end{equation*}
Therefore,
\begin{equation}\label{e-final-7}
\begin{aligned}
&\int_{B_{r_{(k-1)k_0}}\backslash B_{r_{kk_0}}} P_{r_{kk_0}}(x,y)\tilde{g}(y) dy\\
 &\leq (1-\mu)\left(\hat{C}M r_{(k-1)k_0}^{\beta}-M r_{(k-1)k_0}^{\alpha}\right)+M r_{(k-1)k_0}^{\beta}\\
&\leq \hat{C}M r_{(k+1)k_0}^{\beta}\cdot \frac{1-\mu}{\tau_1^{2\beta k_0}}+\mu M r_{(k-1)k_0}^{\beta}\\
&\leq \hat{C}M r_{(k+1)k_0}^{\beta}\left(\frac{1-\mu}{\tau_1^{2\beta k_0}}+\frac{\mu}{\hat{C}\tau_1^{2\beta k_0}}\right),
\end{aligned}
\end{equation}
where $0<\mu<1$ depends only on $n, s,\tau_1,\tau_2$ and $\nu$.

Now, we estimate the second part in\cref{e-final-11} for $x\in B_{r_{(k+1)k_0}}$. Take $k_0$ large enough such that
\begin{equation}\label{e-final-12}
  \tau_2^{k_0}\leq 1/4.
\end{equation}
Then it is easy to check that $|y-\tilde{x}|\leq 2|y-x|$ for any $\tilde{x}\in B_{r_{kk_0}}$ and $y\in B_{r_{(k-1)k_0}}^c$. Let $t=\max(r_{(k+1)k_0}/r_{kk_0},r_{kk_0}/r_{(k-1)k_0})$ and $\tilde{x}=x/t\in B_{r_{kk_0}}$. Then
\begin{equation}\label{e-final-8}
\begin{aligned}
&\int_{B^c_{r_{(k-1)k_0}}} P_{r_{kk_0}}(x,y)\tilde{g}(y) dy=\int_{B^c_{r_{(k-1)k_0}}} P_{r_{kk_0}}(x,y)v_k(y) dy\\
&=C\int_{B^c_{r_{(k-1)k_0}}} \left(\frac{r_{kk_0}^2-|x|^2}{|y|^2-r_{kk_0}^2}\right)^{s/2}\cdot\frac{v_k(y)}{|x-y|^n} dy\\
&=Ct^s\int_{B^c_{r_{(k-1)k_0}}} \left(\frac{r_{kk_0}^2/t^2-|\tilde{x}|^2}{|y|^2-r_{kk_0}^2}\right)^{s/2}\cdot\frac{v_k(y)}{|x-y|^n} dy\\
&\leq Ct^s\int_{B^c_{r_{(k-1)k_0}}} \left(\frac{r_{(k-1)k_0}^2-|\tilde{x}|^2}{|y|^2-r_{(k-1)k_0}^2}\right)^{s/2}\cdot\frac{v_k(y)}{|\tilde{x}-y|^n} dy\\
 &= C_2t^s\int_{B^c_{r_{(k-1)k_0}}} P_{r_{(k-1)k_0}}(\tilde{x},y)v_k(y) dy\\
&\leq C_2t^s \|v_k\|_{L^{\infty}(B_{r_{kk_0}})}\\
&\leq C_2t^s \hat{C}Mr_{kk_0}^{\beta}\\
&= \hat{C}Mr_{(k+1)k_0}^{\beta} \cdot C_2\frac{t^{s}r_{kk_0}^{\beta}}{r_{(k+1)k_0}^{\beta}},
\end{aligned}
\end{equation}
where $C_2$ depends only on $n$. On the other hand,
\begin{equation}\label{e-final-9}
\begin{aligned}
\int_{ B_{r_{kk_0}}} G(x,y)|f|(y)  dy\leq C_1Mr_{kk_0}^{s-n/p}
\leq C_1Mr_{kk_0}^{\beta}
\leq \hat{C}M r_{(k+1)k_0}^{\beta}\cdot\frac{C_1}{\hat{C}\tau_1^{\beta k_0}}.
\end{aligned}
\end{equation}

Take $k_0$ large enough such that\cref{e-final-12} holds and
\begin{equation*}
  C_2\tau_2^{k_0s/2}\leq \mu/8.
\end{equation*}
Then, take $\beta$ small enough such that\cref{e-final-6} holds,
\begin{equation*}
  \frac{1-\mu}{\tau_1^{2\beta k_0}}\leq 1-\mu/2
\end{equation*}
and
\begin{equation*}
  \tau_2^{s/2}\leq \tau_1^{\beta}.
\end{equation*}
Hence,
\begin{equation*}\label{e-final-14}
C_2\frac{t^{s}r_{kk_0}^{\beta}}{r_{(k+1)k_0}^{\beta}} \leq C_2\frac{\tau_2^{k_0s}}{\tau_1^{k_0\beta}}\leq C_2 \tau_2^{k_0s/2}\leq \mu/8.
\end{equation*}
Finally, take $\hat{C}$ large enough such that\cref{e-final-10} holds,
\begin{equation*}\label{e-final-15}
  \frac{\mu}{\hat{C}\tau_1^{2\beta k_0}}\leq \mu/8
\end{equation*}
and
\begin{equation*}\label{e-final-16}
  \frac{C_1}{\hat{C}\tau_1^{\beta k_0}}\leq \mu/8.
\end{equation*}

Thus, combining\crefrange{e-final-v}{e-final-7},\cref{e-final-8} and\cref{e-final-9}, we have
\begin{equation*}
\begin{aligned}
\sup_{ B_{r_{(k+1)k_0}}} v_{k+1} &\leq \hat{C}M r_{(k+1)k_0}^{\beta}.
\end{aligned}
\end{equation*}
Therefore, the\cref{e-final-3} holds for $k+1$. By induction, \cref{e-final-3} holds for any $k\geq 2$ and the proof is completed.~\qed

\section{Solvability of the Dirichlet problem for the Laplace equation}
As it is pointed out above, for the Laplace equation, there exists an explicit relationship between the solution and the boundary value via the Poisson integral. Hence, we have clear quantitative estimate for the decay around some boundary point. This allows us to obtain more general boundary regularity. Here, we derive the continuity of the solution up to the boundary, i.e., solve the Dirichlet problem for the Laplace equation.

Let $g$ be continuous at $x_0$. Then denote its modulus of continuity at $x_0$ by
\begin{equation*}
\omega_{g;x_0}(r)=\sup_{|x-x_0|\leq r} |g(x)-g(x_0)|.
\end{equation*}
Hence, $\omega_{g;x_0}(r)$ is nondecreasing and tends to $0$ as $r\rightarrow 0$. If $x_0$ is the origin, write $\omega_{g}(r)=\omega_{g;0}(r)$ for short.

Now we introduce a modified Perron's method from \cite[Chap. 2.8]{MR1814364}. Let $g$ be a bounded function on $\partial \Omega$ and $v$ be continuous in $\Omega$. The $v$ will be called a subfunction relative to $g$ if for every ball $B\subset R^{n}$ and every harmonic function $h$ in $B$ satisfying
\begin{equation*}
\left\{\begin{aligned}
&\liminf_{x\rightarrow x_0} h(x) \geq v(x_0), &&~~x_0\in \partial B\cap \Omega;\\
&\liminf_{x\rightarrow x_0} h(x) \geq \limsup_{x\rightarrow x_0} g(x), &&~~x_0\in \partial \Omega\cap \bar{B},
\end{aligned}\right.
\end{equation*}
we have $v\leq h$ in $\Omega\cap B$. Define
\begin{equation*}
  S_{g}=\left\{v\in C(\Omega)\big | v \mbox{ is  a subfunction relative to } g\right\}.
\end{equation*}

It is easy to check that if $v\in S_{g}$ and $\bar{v}$ is the harmonic lifting of $v$ in $B$ for some ball $B\subset\subset \Omega$, then $\bar{v}\in S_{g}$. And it can be proved the following (see \cite[
Chap. 2.8]{MR1814364}):
\begin{lemma}\label{l-Lap}
The function $u(x)=\sup_{v\in S_g} v(x)$ is harmonic in $\Omega$.
\end{lemma}

Now, we prove the boundary continuity for the solution of the Laplace equation.
\begin{theorem}\label{t-final}
Suppose that $\Omega$ satisfies the (H4) condition at $0\in \partial \Omega$. Let $g$ be bounded on $\partial \Omega$ and continuous at $0$. Then there exists a harmonic function $u$ in $\Omega$, which is continuous up to $0$ and
\begin{equation*}
  |u(x)-u(0)|\leq \omega(|x|), ~~\forall ~x\in \Omega \cap B_1,
\end{equation*}
where $\omega$ is a modulus of continuity depending only on $\omega_g$ and the quantities in the (H4) condition.
\end{theorem}
\proof Let $u$ be defined as in \Cref{l-Lap}. We only need to prove that $u$ is continuous up to $0$. We assume that $g(0)=0$ as previous and let $\Omega_{r}=\Omega\cap B_{r}$. For $k\geq 0$, let $\Gamma_k=\partial B_{r_k}\cap \Omega^c$, $a_k=H^{n-1}(\Gamma_k)/r_k^{n-1}$,
\begin{equation*}
 \tilde{g}_k(x) \equiv\left\{\begin{aligned}
&0 ~~&&\mbox{on}~~ \Gamma_k;\\
&1 ~~&&\mbox{on}~~\partial B_{r_k}\backslash \Gamma_k,
  \end{aligned}\right.
\end{equation*}
and
\begin{equation*}
  v_k(x)=\frac{r_k^2-|x|^2}{n\omega _n r_k} \int_{\partial B_{r_k}}\frac{\tilde{g}_k(y)ds}{|x-y|^n}.
\end{equation*}

Then
\begin{equation}\label{e-final-2}
  v_k \leq 1-c_0 a_k~~\mbox{on}~~B_{r_{k+1}},
\end{equation}
where $0<c_0<1$ depends only on $n$ and $\tau_2$. By the (H4) condition, $\sum_{k=0}^{\infty}a_k=\infty$, which is equivalent to
\begin{equation*}
  \prod_{k=0}^{\infty}(1-c_0a_k)=0.
\end{equation*}
Let $A_0=\|g\|_{L^{\infty}(\Omega)}$. For $k\geq 1$, define $A_k$ as follows:
\begin{equation*}
A_k=\max\big(\omega_g(r_k), \left(1-c_0a_{k-1}/2\right)A_{k-1}\big).
\end{equation*}
Then it is easy to check that $A_k\rightarrow 0$ decreasingly as $k\rightarrow \infty$.

To prove that $u$ is continuous up to $0$, we only need to prove
\begin{equation}\label{e-final-discrete}
\|u\|_{L^{\infty }(\Omega _{r_{k}})}\leq 2A_k.
\end{equation}
We prove above by induction. For $k=0$,\cref{e-final-discrete} holds clearly. Suppose that it holds for $k$. We need to prove that it holds for $k+1$.

Let $\bar{g}_k=\big(2A_k-\omega_g(r_k)\big)\tilde{g}_k+\omega_g(r_k)$ and
\begin{equation*}
  \bar{v}_k(x)=\frac{r_k^2-|x|^2}{n\omega _n r_k} \int_{\partial B_{r_k}}\frac{\bar{g}_k(y)ds}{|x-y|^n}.
\end{equation*}
Then it is easy to verify as before that
\begin{equation*}
u \leq \bar{v}_k ~~\mbox{in}~~\Omega_{r_k}.
\end{equation*}
On the other hand, it can be checked that $\max (u, -\bar{v}_k)\in S_g$ and hence
\begin{equation*}
  -\bar{v}_k \leq u ~~\mbox{in}~~\Omega_{r_k}.
\end{equation*}

From\cref{e-final-2}, we have for $x\in B_{r_{k+1}}$,
\begin{equation*}
 \bar{v}_k(x) \leq \left(1- c_0a_k\right)\left(2A_k-\omega_g(r_k)\right)+\omega_g(r_k).
\end{equation*}
Hence,
\begin{equation*}
\begin{aligned}
  \|u\|_{L^{\infty}(B_{r_{k+1}})}&\leq \left(1- c_0a_k\right)\left(2A_k-\omega_g(r_k)\right)+\omega_g(r_k)\\
&= 2A_k\left(1- c_0a_k\right)+c_0a_k\omega_g(r_k)\\
&\leq 2A_k\left(1-c_0a_k/2\right)\\
&\leq 2A_{k+1}.
\end{aligned}
\end{equation*}
Thus,\cref{e-final-discrete} holds for $k+1$. By induction, the proof is completed. \qed~\\

An immediate corollary is
\begin{corollary}
Let $\Omega$ satisfy the (H4) condition at every point of $\partial \Omega$. Then the Dirichlet problem
\begin{equation*}
 \left\{ \begin{aligned}
    \Delta u&= 0 ~~\mbox{in}~~\Omega;\\
    u &= g ~~\mbox{on}~~\partial\Omega
  \end{aligned}\right.
\end{equation*}
is uniquely solvable for any continuous function $g$. That is, there exists a unique harmonic function in $\Omega$, which is continuous on $\bar{\Omega}$ and coincides with $g$ on $\partial \Omega$.
\end{corollary}
~\\

\section*{References}
\bibliographystyle{elsarticle-num}

\end{document}